\documentclass{article}
\usepackage{spconf,amsmath,graphicx}

\usepackage[caption=false]{subfig}
\usepackage{amssymb,amsfonts}
\usepackage{textcomp}
\usepackage{xcolor}

\usepackage{here}
\usepackage{enumerate}
\usepackage{bm}
\usepackage{algorithm}
\usepackage{algpseudocode}
\graphicspath{{Fig/}}
\usepackage{amsthm}
\usepackage{comment}
\usepackage{cases}

\usepackage{mathtools}
\usepackage{bm}
\usepackage{comment}
\usepackage{enumitem}
\usepackage{cite}
\usepackage{stmaryrd}
\usepackage{csvsimple}
\usepackage{multicol}
\usepackage{multirow}
\usepackage{array}

\usepackage{mathabx}

\mathtoolsset{showonlyrefs=true}

\newcommand{\argmin}{\mathop{\rm arg~min}\limits}

\newcommand{\inprod}[2]{{\langle #1,#2 \rangle}}
\newcommand{\trace}{{\rm Tr}}

\newcommand{\St}{{\rm St}}
\newcommand{\Gr}{\textrm{Gr}}

\newcommand{\T}{\mathrm{T}}

\newcommand{\prox}[1]{\mathrm{prox}_{#1}}
\newcommand{\moreau}[2]{{}^{#2}#1}

\theoremstyle{plain}

\theoremstyle{definition}
\newtheorem{theorem}{Theorem}[section]
\newtheorem{lemma}[theorem]{Lemma}

\newtheorem{problem}[theorem]{Problem}
\newtheorem{definition}[theorem]{Definition}

\newtheorem{remark}[theorem]{Remark}

\title{A variable smoothing for Nonconvexly constrained nonsmooth optimization with application to sparse spectral clustering}
\name{Keita Kume and Isao Yamada\thanks{This work was supported partially by JSPS Grants-in-Aid (19H04134, 22KJ1270) and by JST SICORP (JPMJSC20C6).}}
\address{Dept. of Information and Communications Engineering, Tokyo Institute of Technology
\\Email:\{kume,isao\}@sp.ce.titech.ac.jp}
\begin{document}
\ninept
\maketitle
\renewcommand{\thefootnote}{\alph{footnote})}
\renewcommand{\appendixname}{Appendix}

\begin{abstract}
  We propose a variable smoothing algorithm for solving nonconvexly constrained nonsmooth optimization problems.
  The target problem has two issues that need to be addressed: (i) the nonconvex constraint and (ii) the nonsmooth term.
  To handle the nonconvex constraint, we translate the target problem into an unconstrained problem by parameterizing the nonconvex constraint in terms of a Euclidean space.
  We show that under a certain condition, these problems are equivalent in view of finding a stationary point.
  To find a stationary point of the parameterized problem, the proposed algorithm performs the gradient descent update for the smoothed version of the parameterized problem with replacement of the nonsmooth function by the Moreau envelope, inspired by a variable smoothing algorithm [Böhm-Wright, J. Optim. Theory Appl., 2021] specialized for unconstrained nonsmooth optimization.
  We also present a convergence analysis of the proposed algorithm as well as its application to a nonconvex reformulation of the sparse spectral clustering.
\end{abstract}
\begin{keywords}
  nonconvex constraint, nonsmooth optimization, weakly convex function, Moreau envelope, sparse spectral clustering
\end{keywords}

\renewcommand{\baselinestretch}{0.82}

\vspace{-0.5em}
\section{Introduction} \label{sec:introduction}
\vspace{-0.5em}
Nonsmooth optimization problems with possibly nonconvex constraint have been serving as key models for finding sparse solutions in wide range of signal processing and machine learning including, e.g., sparse Principal Component Analysis (PCA)~\cite{Benidis-Sun-Babu-Palomar16}, sparse variable PCA~\cite{Ulfarsson-Solo08}, robust sparse PCA~\cite{Breloy-Kumar-Sun-Palomar21}, orthogonal dictionary learning/robust subspace recovering~\cite{Chen-Deng-Ma-So21}, sparse spectral clustering~\cite{Lu-Yan-Lin16}, and robust low-rank matrix completion~\cite{Cambier-Absil16}.
Such applications can be formulated as
\begin{problem}\label{problem:nonsmooth}
  For a locally closed\footnote{
    $C$
    is locally closed if,
    for every
    $\bm{x} \in C$,
    there exists a closed neighborhood
    $\mathcal{N} \subset \mathcal{X}$
    of
    $\bm{x}$
    such that
    $C \cap \mathcal{N}$
    is closed\cite{Rockafellar-Wets98}.
  } (nonconvex) nonempty subset
  $C \subset \mathcal{X}$
  of a Euclidean space
  $\mathcal{X}$,
\begin{equation}
  \mathop{\mathrm{Minimize}}_{\bm{x} \in C}\ f(\bm{x})\coloneqq h(\bm{x}) + g\circ G(\bm{x}), \label{eq:nonsmooth}
\end{equation}
where
(i)
$h:\mathcal{X}\to\mathbb{R}$
is differentiable and
$\nabla h$
is Lipschitz continuous with
$L_{\nabla h}>0$;
(ii)
$g:\mathcal{Z} \to \mathbb{R}$
is (possibly nonsmooth)
$\eta$-weakly convex with
$\eta > 0$,
i.e.,
$g + \frac{\eta}{2}\|\cdot\|^{2}$
is convex, and Lipschitz continuous with
$L_{g} > 0$
over a Euclidean space
$\mathcal{Z}$;
(iii)
$G:\mathcal{X} \to \mathcal{Z}$
is smooth;
(iv)
$\inf_{\bm{x}\in C} f(\bm{x}) >-\infty$.
\end{problem}
\noindent
Despite many applications behind Problem~\ref{problem:nonsmooth}, the problem has two difficulties that need to be addressed: (i) the nonconvex constraint
$C$
and (ii) the nonsmoothness of
$g\circ G$.
To address each of the two difficulties, computational ideas have been proposed respectively.

To deal with the nonconvex constraint
$C \subset \mathcal{X}$
in~\eqref{eq:nonsmooth} with
$g = 0$,
parametrization strategies have been utilized to translate constrained optimization problems~\eqref{eq:nonsmooth}
into unconstrained optimization problems~\cite{Yamada-Ezaki03,Helfrich-Willmott-Ye18,Lezcano19,Kume-Yamada21,Kume-Yamada22,Criscitiello-Boumal22,Levin-Kileel-Boumal22,Kume-Yamada23}.
Indeed, such a constraint set
$C$
appearing in signal processing and machine learning can often be parameterized with a smooth mapping
$F:\mathcal{Y}\to \mathcal{X}$
satisfying
$C=F(\mathcal{Y})$
in terms of a Euclidean space
$\mathcal{Y}$,
e.g., as the so-called {\it trivialization}~\cite{Lezcano19}.
Such parameterizations include, e.g., the Stiefel manifold
$\St(p,N)\coloneqq\{\bm{U} \in \mathbb{R}^{N\times p} \mid \bm{U}^{\T}\bm{U} = \bm{I}\}$~\cite{Yamada-Ezaki03,Kume-Yamada22,Lezcano19},
and the set of all
$m$-by-$n$
low-rank matrices with rank
$r$
at most, to name a few (see~\cite{Levin-Kileel-Boumal22}).
By using
$F$,
we can treat the minimization of
$f:\mathcal{X} \to \mathbb{R}$
over
$C$
as the simpler unconstrained minimization of
$f\circ F:\mathcal{Y} \to \mathbb{R}$.

In a case where
$C = \mathcal{X}$
and
$G = A$
is a linear operator in~\eqref{eq:nonsmooth}, {\it the variable smoothing}~\cite{Bohm-Weight21} considers the minimization of
$h + \moreau{g}{\mu} \circ A$,
instead of
$h + g\circ A$,
with a smoothed surrogate
$\moreau{g}{\mu}:\mathcal{Z}\to \mathbb{R}\ (\mu \in (0,\eta^{-1}))$
of
$g$,
where
$\moreau{g}{\mu}$
is called {\it the Moreau envelope}~\cite{Yamada-Yukawa-Yamagishi11,Bauschke-Combettes17,Abe-Yamagishi-Yamada20,Bauschke-Moursi-Wang21} (see~\eqref{eq:moreau}).
Although there is a gap between
$g$
and
$\moreau{g}{\mu}$,
we have
$\lim_{\mu \to 0} \moreau{g}{\mu}(\bm{z}) = g(\bm{z})\ (\bm{z}\in\mathcal{Z})$.
By controlling the decaying rate of
$\mu$
to
$0$,
the variable smoothing tries to find a minimizer of
$h+g\circ A$
by performing the gradient descent update for minimization of the smooth function
$h+\moreau{g}{\mu}\circ A$.
In a case where
$g$
is {\it prox-friendly}, i.e.,
$\moreau{g}{\mu}$
and
$\nabla \moreau{g}{\mu}$
are available within finite calculations, we do not need to solve any subproblem in the variable smoothing.
Indeed, many functions
$g$
appearing in the real-world applications are prox-friendly, e.g.,
$\ell_{1}$-norm, MCP~\cite{Zhang10}, and SCAD~\cite{Fan-Li01} (see, e.g.,~\cite{prox_repository}).


In this paper, by combining the ideas in the parametrization strategy and the variable smoothing, we resolve the two difficulties in Problem~\ref{problem:nonsmooth} simultaneously:
(i) the nonconvex constraint
$C$;
(ii) the nonsmoothness of
$g\circ G$.
For the first difficulty, we consider a parameterized version of Problem~\ref{problem:nonsmooth} as:
\begin{problem}\label{problem:origin}
  For
  $C$,
  $h$,
  $g$,
  and
  $G$
  defined as in Problem~\ref{problem:nonsmooth},
\begin{equation}
  \mathop{\mathrm{Minimize}}_{\bm{y} \in \mathcal{Y}}\ (h + g\circ G)\circ F(\bm{y})(=f\circ F(\bm{y})), \label{eq:origin}
\end{equation}
where
$F:\mathcal{Y}\to \mathcal{X}$
is a smooth mapping defined over a Euclidean space
$\mathcal{Y}$
such that
$C=F(\mathcal{Y})$.
\end{problem}
\noindent
Theorem~\ref{theorem:necessary_equivalence} in Appendix presents a condition for the equivalence of Problem~\ref{problem:nonsmooth} and Problem~\ref{problem:origin} in view of finding a {\it stationary point} satisfying a necessary condition to be a local minimizer.

For the second difficulty,  we propose to extend the ideas~\cite{Bohm-Weight21} in the variable smoothing.
The proposed algorithm illustrated in Algorithm~\ref{alg:proposed} performs the gradient descent update for the minimization of
$(h+\moreau{g}{\mu_{n}}\circ G) \circ F$
at $n$th update with the Moreau envelope
$\moreau{g}{\mu_{n}}$
of
$g$
with
$\mu_{n} \in (0,\eta^{-1})$
converging to
$0$.
Along with the variable smoothing~\cite{Bohm-Weight21}, the proposed algorithm do not need to solve any subproblem in a case where
$g$
is prox-friendly.
To the authors' best knowledge, the proposed algorithm is the first algorithm to solve the class of Problem~\ref{problem:origin} (see also Remark~\ref{remark:advantage}).
For the proposed algorithm, we present in Theorem~\ref{theorem:convergence_extension} a convergence analysis.
In Section~\ref{sec:numerical}, we present an application with a nonconvex reformulation  of the sparse spectral clustering~\cite{Ng-Michael-Weiss01,Ulrike07,Lu-Yan-Lin16} with numerical experiments.

\noindent\textbf{Notation}:
$\mathbb{N}$
and
$\mathbb{R}$
are respectively the set of all positive integers and of all real numbers.
$\|\cdot\|$
and
$\inprod{\cdot}{\cdot}$
are respectively the Euclidean norm and the standard inner product.
For a smooth mapping
$F:\mathcal{Y} \to \mathcal{X}$,
${\displaystyle \mathrm{D}F(\bm{y}):\mathcal{Y}\to\mathcal{X}:\bm{v}\mapsto \lim_{t\to 0,t\neq 0,t\in\mathbb{R}}\frac{F(\bm{y}+t\bm{v}) - F(\bm{y})}{t}}$
is the G\^{a}teaux derivative of
$F$
at
$\bm{y} \in \mathcal{Y}$,
and it is a linear operator by the smoothness of
$F$.
For a smooth function
$f:\mathcal{X} \to \mathbb{R}$,
$\nabla f$
is the gradient of
$f$,
i.e.,
$\nabla f(\bm{x}) \in \mathcal{X}$
at
$\bm{x} \in \mathcal{X}$
satisfies
$\mathrm{D}f(\bm{x})[\bm{v}] = \inprod{\nabla f(\bm{x})}{\bm{v}}\ (\bm{v}\in\mathcal{X})$
with the G\^{a}teaux derivative of
$f$
at
$\bm{x}$.
For a linear operator
$A:\mathcal{X} \to \mathcal{Z}$,
$A^{*}:\mathcal{Z}\to \mathcal{X}$
is the adjoint of
$A$,
i.e.,
$\inprod{A\bm{x}}{\bm{z}}=\inprod{\bm{x}}{A^{*}\bm{z}}$
holds for all
$\bm{x} \in \mathcal{X}$
and
$\bm{z} \in \mathcal{Z}$.
For
$\bm{y} \in \mathcal{Y}$
and
$S \subset \mathcal{Y}$,
$d(\bm{y},S) := \inf\{\|\bm{y}-\bm{s}\| \in \mathbb{R} \mid \bm{s} \in S\}$
is their distance.

\vspace{-0.5em}
\section{Preliminary}
\vspace{-0.5em}
As extensions of the subdifferential of convex functions, we use the following two subdifferentials of nonconvex functions~\cite{Rockafellar-Wets98}.

\begin{definition}[{\cite[8.3]{Rockafellar-Wets98}}]\label{definition:subdifferential}
  Let
  $f:\mathcal{X} \to \mathbb{R}$
  be a function and
  $\widebar{\bm{x}} \in \mathcal{X}$.
  \begin{enumerate}[label=(\alph*),leftmargin=*,align=left]
    \item
      ${\displaystyle \widehat{\partial} f(\widebar{\bm{x}})\coloneqq \left\{\bm{v} \in \mathcal{X} \mid \liminf_{\bm{x}\to \widebar{\bm{x}},\bm{x}\neq \widebar{\bm{x}}} J_{(\widebar{\bm{x}},\bm{v})}(\bm{x}) \geq 0 \right\}}$\footnote{
        ${\displaystyle\liminf_{\bm{x}\to \widebar{\bm{x}},\bm{x}\neq \widebar{\bm{x}}} J_{(\widebar{\bm{x}},\bm{v})}(\bm{x}) = \inf\left\{a \in \mathbb{R} \mid \substack{\exists (\bm{x}_{n})_{n=1}^{\infty} (\subset \mathcal{X}\setminus \{\widebar{\bm{x}}\})\to \widebar{\bm{x}}\ \mathrm{s.t.}\\ \lim_{n\to\infty}J_{(\widebar{\bm{x}},\bm{v})}(\bm{x}_{n}) = a}\right\}}$.
    }
      with
      $J_{(\widebar{\bm{x}},\bm{v})}(\bm{x})$ $\coloneqq \frac{f(\bm{x}) - f(\widebar{\bm{x}}) - \inprod{\bm{v}}{\bm{x}-\widebar{\bm{x}}}}{\|\bm{x}-\widebar{\bm{x}}\|}$
      is {\it the regular  subdifferential} of
      $f$
      at
      $\widebar{\bm{x}}$.
    \item
      $\partial f(\widebar{\bm{x}})$
      defined below
      is {\it the general subdifferential} of
      $f$
      at
      $\widebar{\bm{x}}$:
      \begin{equation}
        \partial f(\widebar{\bm{x}}) \coloneqq
        \left\{\bm{v} \in \mathcal{X} \mid \substack{\displaystyle\exists (\bm{x}_{n})_{n=1}^{\infty}(\subset \mathcal{X})\to \widebar{\bm{x}}, \exists \bm{v}_{n}\in \widehat{\partial}f(\bm{x}_{n})\\\displaystyle \mathrm{s.t.}\ \lim_{n\to\infty} f(\bm{x}_{n}) = f(\widebar{\bm{x}}), \lim_{n\to\infty} \bm{v}_{n}=\bm{v} }\right\}
      \end{equation}
  \end{enumerate}
\end{definition}

For all functions in this paper,
i.e.,
$h,h\circ F,g,g\circ G,g\circ G \circ F, f, f\circ F$,
we can confirm that the above two subdifferentials coincide by combining facts in~\cite[7.25, 7.27, 9.13, 10.6, 10.9]{Rockafellar-Wets98}.
Since the general subdifferential has good properties for our nonsmooth analysis, we focus on the general subdifferential in the following.

From the nonconvexity of Problem~\ref{problem:origin}, its realistic goal is finding a stationary point satisfying a necessary condition for the optimality.
By Fermat's rule~\cite[10.1]{Rockafellar-Wets98}, the standard optimality condition for
$\bm{y}^{\star} \in \mathcal{Y}$
to be a local minimizer of Problem~\ref{problem:origin} is given as
\begin{align}
  \bm{0} & \in \partial (f\circ F)(\bm{y}^{\star}). \label{eq:optimal_origin}
\end{align}
We present, in Theorem~\ref{theorem:necessary_equivalence} of Appendix, a sufficient condition for the equivalence between necessary conditions for the local optimality of both Problems~\ref{problem:nonsmooth} and~\ref{problem:origin} (Note: the main ideas in this paper can be followed without Appendix).

{\it The proximity operator} and {\it the Moreau envelope} have been used extensively to solve nonsmooth optimization problems~\cite{Yamada-Yukawa-Yamagishi11,Bauschke-Combettes17,Abe-Yamagishi-Yamada20,Bauschke-Moursi-Wang21}.
For an
$\eta$-weakly convex function
$g$
with
$\eta > 0$,
its proximity operator and Moreau envelope with
$\mu \in (0,\eta^{-1})$
are respectively defined as
\begin{align}
  (\widebar{\bm{z}} \in \mathcal{Z}) \quad & \prox{\mu g}(\widebar{\bm{z}})  \coloneqq \argmin_{\bm{z} \in \mathcal{Z}} \left(g(\bm{z}) + \frac{1}{2\mu}\|\bm{z}-\widebar{\bm{z}}\|^{2}\right), \label{eq:prox} \\
  (\widebar{\bm{z}} \in \mathcal{Z}) \quad & \moreau{g}{\mu}(\widebar{\bm{z}}) \coloneqq \inf_{\bm{z} \in \mathcal{Z}} \left(g(\bm{z}) + \frac{1}{2\mu}\|\bm{z}-\widebar{\bm{z}}\|^{2}\right) \leq g(\widebar{\bm{z}}), \label{eq:moreau}
\end{align}
where
$\prox{\mu g}$
is single-valued due to the strong convexity of
$g + (2\mu)^{-1}\|\cdot\|^{2}$.
The Moreau envelope
serves as a smoothed surrogate function of
$g$
because
$\lim_{\mu\to 0}\moreau{g}{\mu}(\bm{z}) = g(\bm{z})\ (\bm{z} \in \mathcal{Z})$,
and
$\moreau{g}{\mu}$
is continuously differentiable with
$\nabla \moreau{g}{\mu}(\bm{z}) = \mu^{-1}(\bm{z}-\prox{\mu g}(\bm{z}))$
(see, e.g.,~\cite{Bauschke-Moursi-Wang21}).
Many functions
$g$
have the close-form expressions of
$\prox{\mu g}$
and
$\moreau{g}{\mu}$,
e.g.,
$\ell_{1}$-norm, MCP, and SCAD (see, e.g.,~\cite{prox_repository}).

\vspace{-0.5em}
\section{Variable smoothing with parametrization}
\vspace{-0.5em}
Since Problem~\ref{problem:origin} is still challenging due to the nonsmoothness and the nonconvexity of
$g\circ G$,
we further reduce Problem~\ref{problem:origin} to the minimization of the smoothed version
$(h+\moreau{g}{\mu} \circ G) \circ F$
of~\eqref{eq:origin}
with the Moreau envelope
$\moreau{g}{\mu}\ (\mu \in (0,\eta^{-1}))$.
The following theorem presents a relation between the optimality conditions of Problem~\ref{problem:origin} in~\eqref{eq:optimal_origin} and its smoothed problem.

\begin{theorem}\label{theorem:stationary}
  Consider Problem~\ref{problem:origin}.
  Let
  $(\mu_{n})_{n=1}^{\infty} \subset \mathbb{R}$
  be a sequence of
  $\mu_{n} \in (0,\eta^{-1})$
  converging to
  $0$,
  and let
  $(\bm{y}_{n})_{n=1}^{\infty} \subset \mathcal{Y}$
  be a sequence converging to some
  $\bm{y}^{\star} \in \mathcal{Y}$.
  For
  $f:=h+g\circ G$,
  we have
  \begin{align} 
     d\left(\bm{0}, \partial (f\circ F)(\bm{y}^{\star})\right)
     \leq \liminf_{n\to\infty}\|\nabla \left((h+\moreau{g}{\mu_{n}}\circ G)\circ F\right)(\bm{y}_{n})\|.\label{eq:necessary_bound}
  \end{align}
\end{theorem}

By noting the condition~\eqref{eq:optimal_origin} for
$\bm{y}^{\star} \in \mathcal{Y}$
can be restated as
\begin{equation}
     d\left(\bm{0}, \partial (f\circ F)(\bm{y}^{\star})\right) = 0,
  \label{eq:optimal_origin_restate}
\end{equation}
Theorem~\ref{theorem:stationary} suggests that a stationary point of Problem~\ref{problem:origin} can be found via the following problem with
$\mu_{n} \in (0,\eta^{-1})$
converging to
$0$
and with
$f_{[n]}\coloneqq h+\moreau{g}{\mu_{n}}\circ G$:
\begin{equation}
  \mathrm{Find}\ 
  (\bm{y}_{n})_{n=1}^{\infty} \subset \mathcal{Y}\ \mathrm{s.t.}\ 
  \liminf_{n\to\infty}\|\nabla (f_{[n]}\circ F)(\bm{y}_{n})\|= 0. \label{eq:stationary}
\end{equation}

\begin{algorithm}[t]
  \caption{Variable smoothing for Problem~\ref{problem:origin}}
  \label{alg:proposed}
\algrenewcommand\algorithmicindent{0.9em}%
  {\footnotesize
    \begin{algorithmic}[0]
      \State \hspace{-1em}
      {\bf Input:}
      $\bm{y}_{1}\in \mathcal{Y},\tau > 2,c\in(0,1),\rho \in (0,1),\alpha > 1, \gamma_{\rm initial} > 0$
    \State \hspace{-1em} {\bf for}
    $n=1,2,\ldots$
    {\bf do}
    \State \hspace{0em}
    $\mu_{n} \leftarrow \frac{1}{\tau\eta\sqrt[\alpha]{n}}$
    \State \hspace{0em}
      Set
      $f_{[n]}\coloneqq h+\moreau{g}{\mu_{n}}\circ G$, and
      $\gamma_{n}\leftarrow \gamma_{\rm initial}$.
    \State \hspace{0em}
        $\bm{y}_{n+1} \leftarrow \bm{y}_{n}-\gamma_{n}\nabla (f_{[n]}\circ F)(\bm{y}_{n})$
        \State \hspace{0em} {\bf while}
        $f_{[n]}(F(\bm{y}_{n+1})) > f_{[n]}(F(\bm{y}_{n})) - c\gamma_{n}\Vert \nabla (f_{[n]}\circ F)(\bm{y}_{n})\Vert^{2}$
        {\bf do}
        \State \hspace{1em}
          $\gamma_{n} \leftarrow \rho \gamma_{n}$
          \Comment{Backtracking algorithm}
        \State \hspace{1em}
          $\bm{y}_{n+1} \leftarrow \bm{y}_{n}-\gamma_{n}\nabla (f_{[n]}\circ F)(\bm{y}_{n})$
      \State \hspace{0em} {\bf end while}
      \State \hspace{-1em}{\bf end for}
  \end{algorithmic}
}
\end{algorithm}

To solve~\eqref{eq:stationary}, we propose to extend the ideas~\cite{Bohm-Weight21} in the variable smoothing.
The proposed algorithm illustrated in Algorithm~\ref{alg:proposed} performs the gradient descent update, from
$\bm{x}_{n} \in \mathcal{X}$
to
$\bm{x}_{n+1} \in \mathcal{X}$,
for the minimization of
$f_{[n]} \circ F= (h+\moreau{g}{\mu_{n}}\circ G) \circ F$,
where the parameter
$\mu_{n} \in (0,\eta^{-1})$
decreases by the rate
$O(n^{-1/\alpha})$
with
$\alpha > 1$.
To find a stepsize
$\gamma_{n} > 0$
reducing the value of the cost function
$f_{[n]} \circ F$
sufficiently, we use the so-called {\it backtracking algorithm} (see, e.g.,~\cite{Nocedal-Wright06}) in Algorithm~\ref{alg:proposed}.

To derive a convergence analysis for Algorithm~\ref{alg:proposed}, we need the Lipschitz continuity of
$\nabla (f_{[n]} \circ F)$.
The following lemma presents its sufficient condition.
\begin{lemma}\label{lemma:extension:Lipschitz_f}
  Consider Problem~\ref{problem:origin}.
  Assume that (i) there exist
  $L_{F},L_{\mathrm{D}F} > 0$
  such that
  $\sup_{\bm{y} \in \mathcal{Y}}\|\mathrm{D}F(\bm{y})\|_{\mathrm{op}} \leq L_{F}$
  and
  $\|\mathrm{D}F(\bm{y}_{1})-\mathrm{D}F(\bm{y}_{2})\|_{\mathrm{op}} \leq L_{\mathrm{D}F}\|\bm{y}_{1} - \bm{y}_{2}\|\ (\bm{y}_{1},\bm{y}_{2} \in \mathcal{Y})$,
  where
  $\|\cdot\|_{\mathrm{op}}$
  is the operator norm;
  (ii) there exist
  $L_{G}, L_{\mathrm{D}G} > 0$
  such that
  $\sup_{\bm{x} \in C}\|\mathrm{D}G(\bm{x})\|_{\mathrm{op}} \leq L_{G}$
  and
  $\|\mathrm{D}G(\bm{x}_{1})-\mathrm{D}G(\bm{x}_{2})\|_{\mathrm{op}} \leq L_{\mathrm{D}G}\|\bm{x}_{1} - \bm{x}_{2}\|\ (\bm{x}_{1},\bm{x}_{2} \in C)$;
  (iii) there exists
  $L_{h}> 0$
  such that
  $\sup_{\bm{x} \in C}\|\nabla h(\bm{x})\| \leq L_{h}$;
  (iv)
  $\mu \in  (0,\min\{1,2^{-1}\eta^{-1}\})$.
  Then,
  $J\coloneqq(h+\moreau{g}{\mu}\circ G)\circ F$
  is continuously differentiable and
  $\nabla J$
  is Lipschitz continuous with
  $\varpi \mu^{-1} > 0$
  over
  $\mathcal{Y}$,
  i.e.,
  \begin{equation} \label{eq:extension:Lipschitz_f_k}
    (\bm{y}_{1},\bm{y}_{2} \in \mathcal{Y}) \quad \|\nabla J(\bm{y}_{1}) - \nabla J(\bm{y}_{2})\| \leq \varpi\mu^{-1}\|\bm{y}_{1}-\bm{y}_{2}\|,
  \end{equation}
  where
$\varpi\coloneqq L_{F}^{2}(L_{\nabla h}+L_{G}^{2}) + L_{\mathrm{D}F}L_{h} + L_{g}(L_{\mathrm{D}F}L_{G}+L_{F}^{2}L_{\mathrm{D}G})$.
\end{lemma}
In a case where
$G$
is a linear operator,
the assumption (ii) on
$G$
in Lemma~\ref{lemma:extension:Lipschitz_f} is satisfied.
Moreover, under the compactness of
$C$,
the assumptions (ii) and (iii) on
$G$
and
$h$
in Lemma~\ref{lemma:extension:Lipschitz_f} are also satisfied automatically.
We will present, in Section~\ref{sec:numerical}, our reformulation of {\it the sparse spectral clustering} satisfies the all assumptions in Lemma~\ref{lemma:extension:Lipschitz_f}.

Consequently, we obtain the following convergence analysis.

\begin{theorem}\label{theorem:convergence_extension}
  Consider Problem~\ref{problem:origin}.
  Let
  $(\bm{y}_{n})_{n=1}^{\infty} \subset \mathcal{Y}$
  be generated by Algorithm~\ref{alg:proposed}.
  Assume that for every
  $n \in \mathbb{N}$,
  $f_{[n]}\circ F\coloneqq(h+\moreau{g}{\mu_{n}}\circ G)\circ F$
  is continuously differentiable and its gradient
  $\nabla (f_{[n]}\circ F)$
  is Lipschitz continuous with
  $L_{\nabla (f_{[n]}\circ F)} \coloneqq \varpi\mu_{n}^{-1} > 0$
  with some common
  $\varpi>0$
  over
  $\mathcal{Y}$.
  Then, if
  $\gamma_{\rm initial}> 2(1-c)L_{\nabla (f_{[1]}\circ F)}^{-1}$,
  we have
  $\liminf_{n\to\infty}\|\nabla (f_{[n]}\circ F)(\bm{y}_{n})\| = 0$.
  More precisely, for some constant
  $\kappa > 0$
  and for all
  $n_{1} > n_{0}$,
  we have
  \begin{equation}
    \min_{n_{0}\leq n \leq n_{1}} \| \nabla (f_{[n]}\circ F)(\bm{y}_{n})\|
    \leq \sqrt{\frac{\kappa }{(n_{1}+1)^{1-\alpha^{-1}}-n_{0}^{1-\alpha^{-1}}}}.
    \label{eq:convergence_rate}
  \end{equation}
  In addition, if a subsequence
  $(\bm{y}_{n(l)})_{l=1}^{\infty}$
  of
  $(\bm{y}_{n})_{n=1}^{\infty}$
  converges to some
  $\bm{y}^{\star} \in \mathcal{Y}$,
  and satisfies
  $\lim_{l\to\infty} \|\nabla (f_{[n(l)]}\circ F)(\bm{y}_{n(l)})\| = 0$,
  then
  $\bm{y}^{\star}$
  is a stationary point of Problem~\ref{problem:origin}, i.e.,
  $\bm{y}^{\star}$
  satisfies~\eqref{eq:optimal_origin}, which is confirmed by Theorem~\ref{theorem:stationary} (Note: the existence of
$\bm{y}^{\star}$
is guaranteed if
$(\bm{y}_{n})_{n=1}^{\infty}$ is bounded\footnote{
$(\bm{y}_{n})_{n=1}^{\infty}$
is bounded if the level set
$\{\bm{y} \in \mathcal{Y} \mid f\circ F(\bm{y}) \leq f_{[1]}\circ F(\bm{y}_{1}) + 2\mu_{1}L_{g}^{2}\}$
is bounded with the initial estimate
$\bm{y}_{1} \in \mathcal{Y}$
in Algorithm~\ref{alg:proposed}.
}).
\end{theorem}

\begin{remark}[Role of $\tau$ in Theorem~\ref{theorem:convergence_extension}]
  A sufficient condition for the Lipschitz assumption of
  $\nabla f_{[n]}\circ F$
  in Theorem~\ref{theorem:convergence_extension} is given in Lemma~\ref{lemma:extension:Lipschitz_f}, where
  $\tau > 2$
  in Algorithm~\ref{alg:proposed} ensures
  $\mu_{n} \in (0, \min\{1,2^{-1}\eta^{-1}\})$
  for the assumption (iv) in Lemma~\ref{lemma:extension:Lipschitz_f}.
\end{remark}

Next remark summarizes advantages of Algorithm~\ref{alg:proposed} over the variable smoothing~\cite{Bohm-Weight21}.
\begin{remark}[Comparisons of~\cite{Bohm-Weight21} and Algorithm~\ref{alg:proposed}]\label{remark:advantage}
  \mbox{}
  \begin{enumerate}[label=(\alph*),leftmargin=*,align=left]
    \item
      The variable smoothing~\cite{Bohm-Weight21} can be applied to only the special instance of Problem~\ref{problem:origin} that
      $C = \mathcal{X} = \mathcal{Y}$,
      $F$
      is the identity operator and
      $G$
      is a linear operator
      $A:\mathcal{X}\to\mathcal{Z}$.
      In contrast, Algorithm~\ref{alg:proposed} can be applied to Problem~\ref{problem:origin} under general setting.
    \item
      The variable smoothing~\cite{Bohm-Weight21} requires the Lipschitz constants
      $L_{\nabla f_{[n]}} > 0$
      of
      $\nabla (h+\moreau{g}{\mu_{n}}\circ A)$
      because every stepsize is given by
      $\gamma_{n}\coloneqq L_{\nabla f_{[n]}}^{-1}$.
      In contrast, Algorithm~\ref{alg:proposed} does not need to know such Lipschitz constants thanks to the backtracking algorithm (see, e.g.,~\cite{Nocedal-Wright06}).
    \item
      For the decaying rate
      $\alpha$
      of
      $\mu_{n}$
      in Algorithm~\ref{alg:proposed}, Algorithm~\ref{alg:proposed} can choose
      $\alpha \in (1,\infty)$
      arbitrarily.
      On the other hand, the corresponding decaying rate in~\cite{Bohm-Weight21} is fixed by
      $3$.
      The difference in the choice of the decaying rate is due to the difference in the goals of convergence analysis.
      The paper~\cite{Bohm-Weight21} aims to derive a better convergence rate with respect to {\it a near stationary point} criterion consisting of a sum of
      $\mu_{n}$
      and
      $\|\nabla (h+\moreau{g}{\mu_{n}}\circ A)(\bm{y}_{n})\|$.
      As
      $\alpha$
      increases, although
      $\|\nabla (h+\moreau{g}{\mu_{n}} \circ A)(\bm{y}_{n})\|$
      decreases rapidly (see~\eqref{eq:convergence_rate}),
      $\mu_{n}$
      decreases slowly.
      To balance the rate of decrease of the two terms,
      $\alpha = 3$
      is used in~\cite{Bohm-Weight21}.
      In contrast, for any
      $\alpha \in (1,\infty)$,
      our convergence analysis in Theorem~\ref{theorem:convergence_extension} ensures to obtain a stationary point directly by finding
      $(\bm{y}_{n})_{n=1}^{\infty}$
      satisfying~\eqref{eq:stationary} (see also~\eqref{eq:necessary_bound} and~\eqref{eq:optimal_origin_restate}).
      According to our experiments, choice of
      $\alpha$
      affects numerical convergence performance of Algorithm~\ref{alg:proposed}, which will be discussed elsewhere.
  \end{enumerate}
\end{remark}

\vspace{-0.5em}
\vspace{-0.5em}
\vspace{-0.5em}

{
\tabcolsep = 0.15em
\begin{table*}[h]
  \vspace{-0.7em}
  \caption{Average scores of NMI~\cite{Strehl-Ghosh03} and ARI~\cite{Hubert-Arabie85}, and their standard deviations (numbers in parentheses).}

  \centering
  {\small
  \begin{tabular}{|c|cc|cc|cc|cc|cc|cc|cc|}\hline
    \multicolumn{1}{|c|}{} & \multicolumn{2}{|c|}{iris} & \multicolumn{2}{|c|}{shuttle} & \multicolumn{2}{|c|}{segmentation} & \multicolumn{2}{|c|}{breast cancer} & \multicolumn{2}{|c|}{glass} & \multicolumn{2}{|c|}{wine} & \multicolumn{2}{|c|}{seeds} \\ \cline{2-15}
      \multicolumn{1}{|c|}{} & \multicolumn{1}{|c}{NMI} &\multicolumn{1}{c|}{ARI}& \multicolumn{1}{|c}{NMI} &\multicolumn{1}{c|}{ARI}& \multicolumn{1}{|c}{NMI} &\multicolumn{1}{c|}{ARI}  & \multicolumn{1}{|c}{NMI} &\multicolumn{1}{c|}{ARI} & \multicolumn{1}{|c}{NMI} &\multicolumn{1}{c|}{ARI} & \multicolumn{1}{|c}{NMI} &\multicolumn{1}{c|}{ARI} & \multicolumn{1}{|c}{NMI} &\multicolumn{1}{c|}{ARI} \\\hline

      \multirow{2}{*}{SC~\cite{Ng-Michael-Weiss01}} & 
      \multicolumn{1}{|c}{0.778} &\multicolumn{1}{c|}{0.745}&
      \multicolumn{1}{|c}{0.387} &\multicolumn{1}{c|}{0.205}&
      \multicolumn{1}{|c}{0.501} &\multicolumn{1}{c|}{0.341}&
      \multicolumn{1}{|c}{0.417} &\multicolumn{1}{c|}{0.419}&
      \multicolumn{1}{|c}{0.321} &\multicolumn{1}{c|}{0.174}&
      \multicolumn{1}{|c}{0.433} &\multicolumn{1}{c|}{0.363}&
      \multicolumn{1}{|c}{0.662} &\multicolumn{1}{c|}{0.659}\\
     & \multicolumn{1}{|c}{(0.000)} &\multicolumn{1}{c|}{(0.000)}&
      \multicolumn{1}{|c}{(0.043)} &\multicolumn{1}{c|}{(0.059)}&
      \multicolumn{1}{|c}{(0.043)} &\multicolumn{1}{c|}{(0.056)}&
      \multicolumn{1}{|c}{(0.000)} &\multicolumn{1}{c|}{(0.000)}&
      \multicolumn{1}{|c}{(0.026)} &\multicolumn{1}{c|}{(0.019)}&
      \multicolumn{1}{|c}{(0.000)} &\multicolumn{1}{c|}{(0.000)}&
      \multicolumn{1}{|c}{(0.007)} &\multicolumn{1}{c|}{(0.010)}  \\\hline

      \multirow{2}{*}{SSC($\ell_{1}$+relax)~\cite{Lu-Yan-Lin16}} & 
      \multicolumn{1}{|c}{0.785} &\multicolumn{1}{c|}{0.786}& 
      \multicolumn{1}{|c}{0.426} &\multicolumn{1}{c|}{0.279}&
      \multicolumn{1}{|c}{0.503} &\multicolumn{1}{c|}{0.343}&
      \multicolumn{1}{|c}{0.433} &\multicolumn{1}{c|}{0.462}&
      \multicolumn{1}{|c}{0.325} &\multicolumn{1}{c|}{0.176}&
      \multicolumn{1}{|c}{0.433} &\multicolumn{1}{c|}{0.363}&
      \multicolumn{1}{|c}{0.671} &\multicolumn{1}{c|}{0.675}\\
     & \multicolumn{1}{|c}{(0.000)} &\multicolumn{1}{c|}{(0.000)}&
      \multicolumn{1}{|c}{(0.047)} &\multicolumn{1}{c|}{(0.063)}&
      \multicolumn{1}{|c}{(0.036)} &\multicolumn{1}{c|}{(0.052)}&
      \multicolumn{1}{|c}{(0.000)} &\multicolumn{1}{c|}{(0.000)}&
      \multicolumn{1}{|c}{(0.021)} &\multicolumn{1}{c|}{(0.016)}&
      \multicolumn{1}{|c}{(0.000)} &\multicolumn{1}{c|}{(0.000)}&
      \multicolumn{1}{|c}{(0.032)} &\multicolumn{1}{c|}{(0.034)}\\\hline

      \multirow{2}{*}{Proposed SSC($\ell_{1}$+Gr)} &
      \multicolumn{1}{|c}{\bf 0.823} &\multicolumn{1}{c|}{\bf 0.818}&
      \multicolumn{1}{|c}{0.427} &\multicolumn{1}{c|}{0.276}&
      \multicolumn{1}{|c}{0.501} &\multicolumn{1}{c|}{0.341}&
      \multicolumn{1}{|c}{0.473} &\multicolumn{1}{c|}{0.547}&
      \multicolumn{1}{|c}{0.323} &\multicolumn{1}{c|}{0.175}&
      \multicolumn{1}{|c}{0.433} &\multicolumn{1}{c|}{0.363}&
      \multicolumn{1}{|c}{0.667} &\multicolumn{1}{c|}{0.668}\\
     & \multicolumn{1}{|c}{(0.000)} &\multicolumn{1}{c|}{(0.000)}&
      \multicolumn{1}{|c}{(0.047)} &\multicolumn{1}{c|}{(0.069)}&
      \multicolumn{1}{|c}{(0.037)} &\multicolumn{1}{c|}{(0.052)}&
      \multicolumn{1}{|c}{(0.000)} &\multicolumn{1}{c|}{(0.000)}&
      \multicolumn{1}{|c}{(0.025)} &\multicolumn{1}{c|}{(0.020)}&
      \multicolumn{1}{|c}{(0.000)} &\multicolumn{1}{c|}{(0.000)}&
      \multicolumn{1}{|c}{(0.000)} &\multicolumn{1}{c|}{(0.000)}\\\hline

      \multirow{2}{*}{Proposed SSC(MCP+Gr)} &
      \multicolumn{1}{|c}{\bf 0.823} &\multicolumn{1}{c|}{\bf 0.818}&
      \multicolumn{1}{|c}{\bf 0.434} &\multicolumn{1}{c|}{\bf 0.294}&
      \multicolumn{1}{|c}{\bf 0.507} &\multicolumn{1}{c|}{\bf 0.351}&
      \multicolumn{1}{|c}{\bf 0.558} &\multicolumn{1}{c|}{\bf 0.664}&
      \multicolumn{1}{|c}{\bf 0.341} &\multicolumn{1}{c|}{\bf 0.180}&
      \multicolumn{1}{|c}{\bf 0.442} &\multicolumn{1}{c|}{\bf 0.376}&
      \multicolumn{1}{|c}{\bf 0.721} &\multicolumn{1}{c|}{\bf 0.708}\\
     & \multicolumn{1}{|c}{(0.000)} &\multicolumn{1}{c|}{(0.000)}&
      \multicolumn{1}{|c}{(0.049)} &\multicolumn{1}{c|}{(0.062)}&
      \multicolumn{1}{|c}{(0.042)} &\multicolumn{1}{c|}{(0.061)}&
      \multicolumn{1}{|c}{(0.000)} &\multicolumn{1}{c|}{(0.000)}&
      \multicolumn{1}{|c}{(0.030)} &\multicolumn{1}{c|}{(0.023)}&
      \multicolumn{1}{|c}{(0.000)} &\multicolumn{1}{c|}{(0.000)}&
      \multicolumn{1}{|c}{(0.037)} &\multicolumn{1}{c|}{(0.044)}\\\hline
    \end{tabular}
    }
  \label{table}
  \vspace{-2em}
\end{table*}
}
\section{Application to sparse spectral clustering}\label{sec:numerical}
\vspace{-0.5em}
\subsection{Reformulation of sparse spectral clustering along with Problem~\ref{problem:origin}}
\vspace{-0.5em}
To evaluate the performance of Algorithm~\ref{alg:proposed}, we carried out numerical experiments in the scenario of {\it the Spectral Clustering (SC)} (see, e.g.,~\cite{Ulrike07}).
The SC is widely used for clustering, e.g., not only clustering on graphs~\cite{Dong-Frossard-Vandergheynst-Nefedov12} in graph signal processing, but also single-cell RNA sequence~\cite{Wang-Liu-Chen-Ma-Xue-Zhao22}, remote sensing image analysis~\cite{Tasdemir-Yalcin-Yildirim15}, and detecting clusters in networks~\cite{Wang-Lin-Wang17}.
To split given data
$(\bm{\xi}_{i})_{i=1}^{N} \subset \mathbb{R}^{d}$
into
$k \leq N$
clusters, the SC is executed in the following steps inspired by the graph theory:
\begin{enumerate}[label=\arabic*.,leftmargin=*,align=left]
  \item
    Construct an affinity matrix
    $\bm{W} \in \mathbb{R}^{N\times N}$
    whose entries
    $\bm{W}_{ij} \geq 0$
    stands for the similarity between
    $\bm{\xi}_{i}$
    and
    $\bm{\xi}_{j}$
    (see, e.g., k-nearest neighborhood graph~\cite{Alshmmari-Stavrakakis-Takatsuka21}).
  \item
    Compute the normalized Laplacian
    $\bm{L} \coloneqq \bm{I} - \bm{D}^{-1/2}\bm{W}\bm{D}^{-1/2} \in \mathbb{R}^{N\times N}$
    by regarding
    $\bm{W}$
    as the adjacency matrix of a certain graph, where
    $\bm{D} \in \mathbb{R}^{N\times N}$
    is the degree (diagonal) matrix of
    $\bm{W}$,
    i.e.,
    $\bm{D}_{ii} = \sum_{j=1}^{N}\bm{W}_{ij}$.
  \item \ 
      \vspace{-0.5em}
    \begin{equation}
      \mathrm{Find}\ \bm{U}^{\star} \in \argmin_{[\bm{U}] \in \Gr(k,N)} \trace(\bm{U}^{\T}\bm{L}\bm{U}), \label{eq:SC}
    \end{equation}
    where
    $\Gr(k,N)\coloneqq\{[\bm{U}] \mid \bm{U} \in \St(k,N)\}$
    is the Grassmann manifold,
    $\St(k,N) \coloneqq \{\bm{U} \in \mathbb{R}^{N\times k} \mid \bm{U}^{\T}\bm{U} = \bm{I}_{k}\}$
    is the Stiefel manifold,
    and
    $[\bm{U}]\coloneqq \{\bm{U}\bm{Q} \in \St(k,N) \mid \bm{Q} \in {\rm O}(k)\coloneqq\St(k,k) \}$
    is the equivalence class with the equivalence relation
    ${\bm{U}\sim\bm{U}'} \overset{\rm def}{\Leftrightarrow} \exists \bm{Q} \in {\rm O}(k)\ \mathrm{s.t.}\ \bm{U}' = \bm{U}\bm{Q}$.\footnote{
    $\trace(\bm{U}^{\T}\bm{L}\bm{U}) = \trace(\bm{U}'^{\T}\bm{L}\bm{U}')$
    holds for
    $\bm{U} \in \St(k,N)$
    and
    $\bm{U}' \in [\bm{U}]$.
  }
  \item
    Form
    $\widehat{\bm{U}}^{\star} \in \mathbb{R}^{N\times k}$
    by normalizing each row of
    $\bm{U}^{\star}$
    to length
    $1$.
  \item
    Treat each row
    $\widehat{\bm{u}}^{\star}_{i} \in \mathbb{R}^{k}$
    of
    $\widehat{\bm{U}}^{\star}$
    as a feature vector of
    $\bm{\xi}_{i}$,
    and cluster
    $(\widehat{\bm{u}}_{i})_{i=1}^{N}$
    into
    $k$
    clusters by the k-means algorithm.
\end{enumerate}
$\bm{U}^{\star}$
in~\eqref{eq:SC} can be obtained by the eigenvalue decomposition of
$\bm{L}$,
and its column vectors are eigenvectors of
$\bm{L}$
corresponding to its $k$th smallest eigenvalues.
As suggested in~\cite{Lu-Yan-Lin16}, in the ideal case where
$\bm{W}$
is block diagonal, i.e.,
$\bm{W}_{ij} = 0$
if
$\bm{\xi}_{i}$
and
$\bm{\xi}_{j}$
belong to different clusters,
$\bm{U}^{\star}\bm{U}^{\star\T} \in \mathbb{R}^{N\times N}$
is block diagonal.
To capture the block diagonal structure of
$\bm{U}^{\star}\bm{U}^{\star\T}$,
the sparsity of
$\bm{U}^{\star}\bm{U}^{\star\T}$
can be naturally assumed because its block off-diagonal matrices are desired to be zero.
This motivates {\it the Sparse SC (SSC)}~\cite{Lu-Yan-Lin16},
where the problem in~\eqref{eq:SC} is replaced by the following problem with
$\lambda > 0$:
\begin{equation}
  \mathrm{Find}\ \bm{U}^{\star} \in \argmin_{[\bm{U}] \in \Gr(k,N)} \trace(\bm{U}^{\T}\bm{L}\bm{U}) + \lambda r(\bm{U}\bm{U}^{\T}), \label{eq:SSC}
\end{equation}
where
$r:\mathbb{R}^{N\times N} \to \mathbb{R}$
is a regularizer to promote the sparsity of solutions.
With
$\ell_{1}$-norm
$r=\|\cdot\|_{1}$,
the SSC in~\cite{Lu-Yan-Lin16} solves a convex relaxation of~\eqref{eq:SSC},
where
$\Gr(k,N)$
is relaxed to a certain closed convex set.

The problem in~\eqref{eq:SSC} is a special instance of Problem~\ref{problem:nonsmooth} if
$r$
is weakly convex and Lipschitz continuous.
Fortunately,
$\Gr(k,N)$
has a smooth parametrization called
{\it the generalized left-localized Cayley transform\footnote{
  $\Psi_{\bm{S}}$
  was proposed originally for parametrization of
  $\St(p,N)$,
  and there exist so-called {\it singular-points} in
  $\St(p,N)$,
  where we can not parameterize by
  $\Psi_{\bm{S}}$.
  However, for parametrization of
  $\Gr(k,N)$,
  we can parameterize all points in
  $\Gr(k,N)$
  by
  $\Psi_{\bm{S}}$ (see~\cite[Theorem 3.2]{Kume-Yamada22}).
}}
$\Psi_{\bm{S}}$
with
$\bm{S} \in {\rm O}(N) =\St(N,N)$~\cite{Kume-Yamada22}:
\begin{equation}
  \Psi_{\bm{S}}:Q_{N,k} \to \Gr(k,N): \bm{V} \mapsto [\bm{S}(\bm{I}-\bm{V})(\bm{I}+\bm{V})^{-1}\bm{I}_{N\times k}], \label{eq:Cayley_inv}
\end{equation}
where
$\bm{I}_{N\times k} \in \mathbb{R}^{N\times k}$
is the first $k$th left block matrix of
$\bm{I}$
and
${\displaystyle Q_{N,k}\coloneqq\left\{\begin{bmatrix} \bm{A} & -\bm{B}^{\T} \\ \bm{B} & \bm{0} \end{bmatrix} \in \mathbb{R}^{N\times N} \mid \substack{\bm{A}^{\T} = -\bm{A} \in \mathbb{R}^{k\times k},\\ \bm{B} \in \mathbb{R}^{(N-k)\times k}}\right\}}$
is a Euclidean space.
With
$h_{\rm SSC}(\bm{U})\coloneqq\trace(\bm{U}^{\T}\bm{L}\bm{U})$
and
$G_{\rm SSC}(\bm{U})\coloneqq\bm{U}\bm{U}^{\T}$,
we can reformulate~\eqref{eq:SSC} along with Problem~\ref{problem:origin}
as
\begin{equation}
  \mathrm{Find}\ \bm{V}^{\star} \in \argmin_{\bm{V} \in Q_{N,k}} (h_{\rm SSC} + \lambda r \circ G_{\rm SSC})\circ \Psi_{\bm{S}}(\bm{V}). \label{eq:SSC_parametrizaiton} 
\end{equation}
\vspace{-0.5em}

All assumptions on functions and mappings in Lemma~\ref{lemma:extension:Lipschitz_f} are satisfied with
$h=h_{\rm SSC}$,
$G=G_{\rm SSC}$,
and
$F=\Psi_{\bm{S}}$,
where the assumption (i) can be confirmed by using~\cite[Proposition A.11, Lemma A.12]{Kume-Yamada22},
the assumptions (ii) and (iii) can be confirmed by the compactness of
$\Gr(k,N)$.
By Theorem~\ref{theorem:convergence_extension}, Algorithm~\ref{alg:proposed} can solve~\eqref{eq:SSC}, via~\eqref{eq:SSC_parametrizaiton}\footnote{
  We can find a stationary point of~\eqref{eq:SSC} by finding a stationary point of~\eqref{eq:SSC_parametrizaiton} since
  $C\coloneqq\Gr(k,N)$
  and
  $F\coloneqq\Psi_{\bm{S}}$
  satisfy the assumptions in Theorem~\ref{theorem:necessary_equivalence}.
}, without any convex relaxation of
$\Gr(k,N)$,
which was used in~\cite{Lu-Yan-Lin16}.
In particular, Algorithm~\ref{alg:proposed} can solve the problem~\eqref{eq:SSC} with weakly convex penalties
$r$
(Note: \cite{Wang-Liu-Chen-Ma-Xue-Zhao22} reported a manifold proximal linear method for solving the problem~\eqref{eq:SSC} in a case where
$r$
is $\ell_{1}$-norm).
Indeed, (Lipschitz continuous) weakly convex penalties, e.g.,
MCP~\cite{Zhang10} and SCAD~\cite{Fan-Li01} (see also~\cite{Bohm-Weight21}), have been utilized to promote more sparsity for solutions than
$\ell_{1}$-norm~\cite{Abe-Yamagishi-Yamada20}.
In the next subsection, we demonstrate the efficacies of the model~\eqref{eq:SSC_parametrizaiton} with a weakly convex penalty and of Algorithm~\ref{alg:proposed}.

\vspace{-1em}
\subsection{Numerical experiments}
\vspace{-0.5em}
We compared the SC~\cite{Ng-Michael-Weiss01}, the SSC with the convex relaxation~\cite{Lu-Yan-Lin16,SSC_convex_github} denoted by SSC($\ell_{1}$+relax), and the SSC by the reformulation~\eqref{eq:SSC_parametrizaiton} with
$\ell_{1}$-norm and MCP~\cite{Zhang10}, denoted respectively by SSC($\ell_{1}$+Gr) and by SSC(MCP+Gr).
More precisely,
as
$r$
in~\eqref{eq:SSC_parametrizaiton}, we used
$\ell_{1}$-norm:
$\|\cdot\|_{1}:=\sum_{i}|z_{i}|$,
and MCP\footnote{
  $g_{{\rm MCP}(\beta)}:\mathbb{R}\to\mathbb{R}:z \mapsto
\begin{cases}
  |z| - \frac{z^{2}}{2\beta} & |z|\leq \beta \\
  \frac{\beta}{2} & |z| > \beta
\end{cases}$
}:
$\sum_{i} g_{{\rm MCP}(\beta)}(z_{i})$
with
$\beta > 0$.

For SSC($\ell_{1}$+relax), SSC($\ell_{1}$+Gr) and SSC(MCP+Gr), parameters were searched from
$\lambda, \beta \in \{10^{-i}\mid i = 0,1,2,\ldots,6\}$.
For SSC($\ell_{1}$+Gr) and SSC(MCP+Gr), we employed Algorithm~\ref{alg:proposed} with parameters
$(\tau,c,\rho,\alpha) = (2.1,2^{-13},0.5,2)$,
and
$\gamma_{\rm initial}:= \max\{1,\Vert \nabla (f_{[1]}\circ F)(\bm{y}_{1})\Vert^{-1}\}$.
We used
$\bm{S} \in {\rm O}(N)$
for
$\Psi_{\bm{S}}$
by using a selection suggested in~\cite[Alg. 1]{Kume-Yamada23}.
We implemented the code in MATLAB (R2023b) on a MacBook Pro (13-inch, M1, 2020) with 16GB of RAM.
We stopped
Algorithm~\ref{alg:proposed} was terminated at
$n$th iterate achieving
$n=2000$
or
$\vert f\circ F(\bm{y}_{n})-f\circ F(\bm{y}_{n-1})\vert < 10^{-8}$.

We evaluated performances of algorithms by using 7 real-world datasets in UCI Machine Learning Repository:
(i) "iris"; (ii) "shuttle" (chosen $1500$ samples randomly); (iii) "segmentation";
(iv) "breast cancer";
(v) "glass";
(vi) "wine";
and (vii) "seeds".
For each dataset, we made an affinity matrix
$\bm{W}$
by following~\cite{Alshmmari-Stavrakakis-Takatsuka21}.
To evaluate the performance, we computed the Normalized Mutual Information (NMI)~\cite{Strehl-Ghosh03} and the Adjusted Rand Index (ARI)~\cite{Hubert-Arabie85}.
These higher scores indicate better clustering performance.

Table~\ref{table} shows the averaged results of each algorithm by running k-means algorithm for $100$ times in the step 5 of the spectral clustering.
From Table~\ref{table}, we can see that the proposed SSC($\ell_{1}$+Gr) is comparable or superior to SSC($\ell_{1}$+relax), implying thus solving~\eqref{eq:SSC} without any relaxation improves the performance of clustering.
Moreover, the proposed SSC(MCP+Gr) has the best performance for every dataset.
This implies that weakly convex penalties such as MCP can achieve further improvements than the convex penalty $\ell_{1}$-norm.
Consequently, the result shows the great potential of Algorithm~\ref{alg:proposed}, and of the weakly convex model~\eqref{eq:SSC_parametrizaiton} along with Problem~\ref{problem:origin}.

\vspace{-1em}
\section*{\normalsize Appendix: Optimality conditions for Problems~\ref{problem:nonsmooth} and~\ref{problem:origin}}
\vspace{-0.5em}
{
In this appendix, we present a relation regarding stationary points between Problem~\ref{problem:nonsmooth} and Problem~\ref{problem:origin} in more general cases where
$f$
belongs to
$\mathfrak{F}(\mathcal{X})$,
where
$\mathfrak{F}(\mathcal{X})$
is the set of all locally Lipschitz continuous functions\footnote{
  $\widetilde{f}:\mathcal{X}\to \mathbb{R}$
  is locally Lipschitz continuous if, for all
  $\widebar{\bm{x}} \in \mathcal{X}$,
  there exists a neighborhood
  $\mathcal{N} \subset \mathcal{X}$
  at
  $\widebar{\bm{x}}$
  such that
  $\widetilde{f}$
  is Lipschitz continuous over
  $\mathcal{N}$.
}
$\widetilde{f}:\mathcal{X} \to \mathbb{R}$
satisfying
$\partial  \widetilde{f}(\widebar{\bm{x}}) = \widehat{\partial} \widetilde{f}(\widebar{\bm{x}})\ (\widebar{\bm{x}}\in \mathcal{X})$.
Clearly
$f:=h+g\circ G$
in Problem~\ref{problem:nonsmooth} belongs to
$\mathfrak{F}(\mathcal{X})$
because
$h$
is continuously differentiable, and
$g$
is Lipschitz continuous.

Under {\it the Clarke regularity}\footnote{ \label{foot:Clarke}
  $C \subset \mathcal{X}$
  is Clarke regular if
  $\widehat{N}_{C}(\widebar{\bm{x}}) = N_{C}(\widebar{\bm{x}})$
  for all
  $\widebar{\bm{x}} \in C$~\cite[6.4]{Rockafellar-Wets98},
  where
  $N_{C}(\bm{x}^{\star})\coloneqq \partial \iota_{C}(\bm{x}^{\star})$
  and
  $\widehat{N}_{C}(\widebar{\bm{x}})\coloneqq\widehat{\partial}\iota_{C}(\widebar{\bm{x}})$
  are {\it the general normal cone} and {\it the regular normal cone} to
  $C$
  at
  $\widebar{\bm{x}} \in C$~\cite[8.14]{Rockafellar-Wets98} respectively,
  and
  $\iota_{C}(\widebar{\bm{x}})\coloneqq\begin{cases} 0 & (\widebar{\bm{x}}\in C) \\ \infty & (\widebar{\bm{x}}\notin C)\end{cases}$
  is {\it the indicator function} of
  $C$.
} of
$C\subset \mathcal{X}$,
a necessary condition for
$\bm{x}^{\star} \in \mathcal{X}$
to be a local minimizer of Problem~\ref{problem:nonsmooth} is given as
\begin{equation}
  \bm{0} \in \partial f(\bm{x}^{\star}) + N_{C}(\bm{x}^{\star}) \label{eq:stationary_nonsmooth}
\end{equation}
by Fermat's rule~\cite[10.1]{Rockafellar-Wets98}, where
$N_{C}(\bm{x}^{\star})$
is {\it the general normal cone}\footref{foot:Clarke} to
$C$
at
$\bm{x}^{\star}$\cite[8.14]{Rockafellar-Wets98}.
Then, we have the following relations between two conditions~\eqref{eq:stationary_nonsmooth} and~\eqref{eq:optimal_origin} for stationary points of Problem~\ref{problem:nonsmooth} and Problem~\ref{problem:origin}.
\begin{theorem}\label{theorem:necessary_equivalence}
  Let
  $\mathcal{X},\mathcal{Y}$
  be Euclidean spaces, and
  $F:\mathcal{Y}\to \mathcal{X}$
  smooth such that
  $(\emptyset \neq)F(\mathcal{Y}) = C \subset \mathcal{X}$.
  Suppose that
  $C$
  is Clarke regular.
  Then, the following hold:
  \begin{enumerate}[label=(\alph*)]
    \item
      For
      $f\in \mathfrak{F}(\mathcal{X})$
      (especially in Problem~\ref{problem:nonsmooth}),
      we have
      \begin{equation}
        \hspace{-1.2em}
        (\forall\bm{y}\in\mathcal{Y})\ 
        \underbrace{\bm{0} \in \partial f(F(\bm{y})) + N_{C}(F(\bm{y}))}_{\substack{\rm Optimality\ condition\\\rm for\ Problem~\ref{problem:nonsmooth}  (see~\eqref{eq:stationary_nonsmooth})}}
        \hspace{-0.2em}\Rightarrow \hspace{-0.7em}\underbrace{\bm{0} \in \partial (f\circ F)(\bm{y}).}_{ \substack{\rm Optimality\ condition\\\rm for\ Problem~\ref{problem:origin}\ (see~\eqref{eq:optimal_origin})}} \label{eq:optimal_equivalence}
  \vspace{-0.4em}
      \end{equation}
    \item The following conditions are equivalent:
    \begin{enumerate}[label=(\roman*)]
      \item
        $(\forall f \in \mathfrak{F}(\mathcal{X}), \forall \bm{y} \in \mathcal{Y})$
        \begin{equation}
          \quad \Bigl(\bm{0} \in \partial f(F(\bm{y})) + N_{C}(F(\bm{y})) \Leftrightarrow  \bm{0} \in \partial (f\circ F)(\bm{y}) \Bigr).
        \end{equation}
      \item
        $(\forall \bm{y} \in \mathcal{Y})\ N_{C}(F(\bm{y})) =  \mathrm{Ker}( (\mathrm{D}F(\bm{y}))^{*})$,
        where
        $\mathrm{Ker}(\cdot)$
        stands for the null space of a linear operator.
    \end{enumerate}
  \end{enumerate}
\end{theorem}

From Theorem~\ref{theorem:necessary_equivalence}, for
$f$
in Problem~\ref{problem:nonsmooth}, if we have
$N_{C}(F(\bm{y}^{\star})) =  \mathrm{Ker}( (\mathrm{D}F(\bm{y}^{\star}))^{*})$
for
$\bm{y}^{\star} \in \mathcal{Y}$,
then the conditions~\eqref{eq:stationary_nonsmooth} and~\eqref{eq:optimal_origin} are equivalent under
$\bm{x}^{\star} = F(\bm{y}^{\star})$.
The condition
$N_{C}(F(\bm{y}^{\star})) =  \mathrm{Ker}( (\mathrm{D}F(\bm{y}^{\star}))^{*})$
is guaranteed in a case where
$C$
is an {\it embedded submanifold} in
$\mathcal{X}$
and
$F$
is {\it submersion}\footnote{
$C$
is a $k$-dimensional embedded submanifold in the
Euclidean space
$\mathcal{X}$
if (i) for each
$\widebar{\bm{x}} \in C$,
there exists an open neighborhood
$\mathcal{N}$
of
$\widebar{\bm{x}}$
in
$\mathcal{X}$
and a smooth mapping
$H:\mathcal{N} \to \mathbb{R}^{k}$
such that
$C\cap \mathcal{N} = H^{-1}(\bm{0})$
and (ii) the rank of
$\mathrm{D}H(\widebar{\bm{x}})$
is
$k$~\cite[Def. 3.10]{Boumal20}.
Moreover,
$F:\mathcal{Y}\to \mathcal{X}$
satisfying
$F(\mathcal{Y}) = C$
is submersion if the rank of
$\mathrm{D}F(\bm{y})$
is
$k$
for
$\bm{y} \in \mathcal{Y}$~\cite[p.211]{Boumal20}.
}.
Finally, we remark that Theorem~\ref{theorem:necessary_equivalence} (b) serves as an extension of a characterization in~\cite[Theorem 2.4]{Levin-Kileel-Boumal22} of stationary points for a smooth function to that for a nonsmooth function.
}

\newpage

\bibliographystyle{IEEEbib}
\bibliography{refs}

\end{document}